\documentclass{amsart}
\usepackage{amssymb}
\usepackage{a4}
\begin{document}
\newtheorem{theorem}{Theorem}[section]
\newtheorem{lemma}[theorem]{Lemma}
\newtheorem{remark}[theorem]{Remark}
\newtheorem{definition}[theorem]{Definition}
\newtheorem{corollary}[theorem]{Corollary}
\newtheorem{example}[theorem]{Example}
\newtheorem{conjecture}[theorem]{Conjecture}
\font\pbglie=eufm10
\def\Tr{\operatorname{Tr}}
\def\id{\operatorname{Id}}
\def\ffrac#1#2{{\textstyle\frac{#1}{#2}}}

\makeatletter
 \renewcommand{\theequation}{%
 \thesection.\alph{equation}}
 \@addtoreset{equation}{section}
 \makeatother
\title[Conformally Osserman manifolds]
{Conformally Osserman manifolds and conformally complex space forms}
\author{N. Bla{\v z}i{\'c} and  P. Gilkey}
\begin{address}{NB: Faculty of Mathematics, University of Beograd, Studentski Trg. 16, P.P. 550,
11000 Beograd, Srbija i Crna Gora. Email: {\it blazicn@matf.bg.ac.yu}}\end{address}
\begin{address}{PG: Mathematics Department, University of Oregon,
Eugene Or 97403 USA.\newline Email: {\it gilkey@darkwing.uoregon.edu}}
\end{address}

\begin{abstract} We characterize manifolds which are locally conformally equivalent to either complex projective
space or to its negative curvature dual in terms of their Weyl curvature tensor.
As a byproduct of this investigation, we classify the conformally complex space forms if the dimension is at least 8. We
also study when the Jacobi operator associated to the Weyl conformal curvature tensor of a Riemannian
manifold has constant eigenvalues on the bundle of unit  tangent
vectors and classify such manifolds which are not conformally flat in dimensions congruent to 2 mod 4.
\end{abstract}
\keywords{complex space form, conformally flat, conformally complex space form, conformally Osserman manifold,
conformal Jacobi operator,  Jacobi operator, Osserman manifold, space form, Weyl conformal tensor.
\newline \phantom{.....}2000 {\it Mathematics Subject Classification.} 53B20.}
\maketitle

\section{Introduction}\label{Sect-1}

\subsection{The Weyl curvature} Let $\nabla$ be the Levi-Civita connection of a Riemannian manifold
$(M,g)$ of dimension
$m$. The curvature operator $R(x,y)$ and curvature tensor $R(x,y,z,w)$ are defined by setting:
$$R(x,y):=\nabla_x\nabla_y-\nabla_y\nabla_x-\nabla_{[x,y]}\quad\text{and}\quad
R(x,y,z,w)=g(R(x,y)z,w)\,.
$$
Let
$\{e_i\}$ be a local orthonormal frame for the tangent bundle. We sum over repeated indices to define the {\it
Ricci tensor} $\rho$ and the {\it scalar curvature} $\tau$ by setting:
$$\rho_{ij}:=\textstyle\sum_kR_{ikkj}\quad\text{and}\quad\tau:=\textstyle\sum_i\rho_{ii}\,.$$
The associated {\it Ricci operator} is defined by setting $\rho(e_i)=\textstyle\sum_j\rho_{ij}e_j$.
We introduce additional tensors by setting
\begin{equation}\label{eqn-1.a}
\begin{array}{l}
L(x,y)z:=g(\rho y,z)x-g(\rho x,z)y+g(y,z)\rho x-g(x,z)\rho y,\\
R_0(x,y)z:=g(y,z)x-g(x,z)y\,.\vphantom{\vrule height 11pt}
\end{array}\end{equation}
Let $W$ be the {\it Weyl conformal curvature operator}. We may decompose
\begin{equation}\label{eqn-1.b}\begin{array}{l}
R(x,y)=W(x,y)+c_1(m)\tau R_0(x,y)+c_2(m)L(x,y)\quad\text{where}\\
c_1(m):=-\ffrac1{(m-1)(m-2)}\quad\text{and}\quad c_2(m):=\ffrac1{m-2}\,.\vphantom{\vrule height 12pt}
\end{array}\end{equation}

\subsection{Conformal geometry} We say that two Riemannian metrics $g_1$ and $g_2$ are {\it conformally
equivalent} if $g_1=\alpha\cdot g_2$ where $\alpha$ is a smooth positive scaling function. The
Weyl conformal curvature operator is invariant on a conformal class as
\begin{equation}\label{eqn-1.c}
W_{g_1}=W_{g_2}\quad\text{if}\quad g_1=\alpha\cdot g_2\,.
\end{equation}

Conformal analogues of notions in Riemannian geometry can be obtained by replacing the full curvature operator
$R$ by the Weyl operator $W$; we add the prefix ``conformally'' in doing this. For example, one says that
$(M,g)$ is {\it conformally flat} if the Weyl tensor $W$ vanishes identically; this implies that $(M,g)$ is
conformally equivalent to flat space.

\subsection{Space forms and complex space forms} One says that $(M,g)$ is a {\it space form} if $R=\lambda_0R_0$ for some smooth function
$\lambda_0$ or, equivalently, if $(M,g)$ has pointwise constant sectional curvature. If $m\ge3$, then necessarily
$\lambda_0$ is constant and by rescaling the metric, we may assume
$\lambda_0\in\{-1,0,1\}$. If $\lambda_0=-1$, then $(M,g)$ is locally isometric to hyperbolic space; if
$\lambda_0=0$, then $(M,g)$ is locally isometric to flat space; if $\lambda_0=1$, then $(M,g)$ is locally
isometric to the sphere. Thus the geometry is very rigid in this setting.

Let $\Phi$ be a Hermitian almost complex structure on $TM$; necessarily $m=2n$ is even. We set
\begin{equation}\label{eqn-1.d}
R_\Phi(x,y)z:=g(\Phi y,z)\Phi x-g(\Phi x,z)\Phi y-2g(\Phi x,y)\Phi z\,.
\end{equation}
We say that $(M,g)$ is a {\it complex space form} if $R=\lambda_0R_0+\lambda_1R_\Phi$ for smooth functions
$\lambda_0$ and $\lambda_1$ where $\lambda_1\ne0$.
Let $(\mathbb{CP}^n,g_{FS})$ denote complex projective
space with the Fubini-Study metric and let $({}^*\mathbb{CP}^n,{}^*g_{FS})$ be the negative curvature
dual; these are complex space forms. Conversely, if $(M,g)$ is a complex space form and if
$m\ge6$, then one can show that
$\lambda_0=\lambda_1$ and that $\lambda_0$  is constant. By rescaling the metric, we may assume $\lambda_0=\pm1$.
If $\lambda_0=1$, then
$(M,g)$ is locally isometric to $(\mathbb{CP}^n,g_{FS})$; if $\lambda_0=-1$, then $(M,g)$ is locally isometric to
$({}^*\mathbb{CP}^n,{}^*g_{FS})$. We refer to \cite{TV81} for further details. A generalization of
this result to the pseudo-Riemannian setting may be found in \cite{Gi02}.

\subsection{Conformally complex space forms} One says that $(M,g)$ is a {\it conformal
space form} if $W=\lambda_0R_0$; we shall see presently this implies $\lambda_0=0$ so $(M,g)$ is conformally flat.
Similarly one says that $(M,g)$ is a {\it conformally complex space form} if
$W=\lambda_0R_0+\lambda_1R_\Phi$ for some Hermitian almost complex structure on $TM$ where $\lambda_0$ and
$\lambda_1$ are smooth functions on $M$ with $\lambda_1\ne0$. A similar rigidity
result holds:

\begin{theorem}\label{thm-1.1}
 Let $(M,g)$ be a conformally complex space form with $m\ge8$. Then $(M,g)$ is
locally conformally equivalent to either $(\mathbb{CP}^n,g_{FS})$ or $({}^*\mathbb{CP}^n,{}^*g_{FS})$.
\end{theorem}

We remark that our proof of Theorem \ref{thm-1.1} extends to the higher signature setting; we shall omit details
in the interests of brevity.

\subsection{The Jacobi operator $J_R$} We define a self-adjoint map $J_R(x)$ of the tangent bundle
defined by setting:
$$J_R(x)y=R(y,x)x\,.$$
One says that $(M,g)$ is {\it Osserman} if the eigenvalues of $J_R(x)$ are
constant on the sphere bundle $S(M,g)$ of unit tangent vectors. One says that $(M,g)$ is a
{\it local $2$ point homogeneous space} if the local isometries of $(M,g)$ act transitively on
$S(M,g)$; this necessarily implies that $(M,g)$ is Osserman. Osserman \cite{refOss} wondered if the
converse held; this has been called the Osserman conjecture by subsequent
authors. Chi
\cite{Ch88} and Nikolayevsky
\cite{refNik,refNik2} established the Osserman conjecture if $m\ne8,16$.

There are similar questions for pseudo-Riemannian manifolds.  Let $(M,g)$ be a pseudo-Riemannian manifold of
signature $(p,q)$. One says that $(M,g)$ is spacelike (resp. timelike) Jordan Osserman if the Jordan normal form
of $J_R$ is constant on the pseudo-sphere of unit spacelike (resp. timelike) vectors. In the Lorentzian setting
($q=1$), it is known that any spacelike (resp. timelike) Jordan Osserman manifold has constant sectional curvature
\cite{refBBG,refGKV}. The analogous question in the higher signature setting is far from settled. For example,
there are spacelike and timelike Jordan Osserman pseudo-Riemannian manifolds which are not locally homogeneous and
thus are not local symmetric spaces \cite{GKV02}.

\subsection{The conformal Jacobi operator} We follow
the discussion of
\cite{BGNS03} and study a conformal analogue of the Osserman conjecture.
The {\it conformal Jacobi operator} $J_W$ is given by:
$$J_W(x)y:=W(y,x)x;$$
the constants $c_i$ are chosen so that
\begin{equation}\label{eqn-1.e}
\Tr\{J_W\}=0\,.
\end{equation}
Thus, in particular, if $W=\lambda_0R_0$, then $0=\Tr(J_W(x))=(m-1)\lambda_0g(x,x)$ so $\lambda_0=0$ and $(M,g)$
is conformally flat as noted above.

We say that
$(M,g)$ is {\it conformally Osserman} if the eigenvalues of
$J_W$ are constant on the fiber spheres
$$S_P(M,g):=\{x\in T_PM:g(x,x)=1\};$$
the eigenvalues are allowed to vary from point to point.
Since space forms and complex space forms are $2$ point homogeneous spaces, they are examples of conformally Osserman
manifolds. In particular, $(\mathbb{CP}^n,g_{FS})$ and $({}^*\mathbb{CP}^n,{}^*g_{FS})$ are conformally Osserman
manifolds.

If $g_1=\alpha g_2$, then
$S_P(M,g_1)=\alpha^{-1/2}S_P(M,g_2)$. Furthermore
$J_{W(g_1)}=J_{W(g_2)}$ by Equation (\ref{eqn-1.c}).
Thus the eigenvalues on the unit sphere bundles rescale so:
\begin{theorem}\label{thm-1.2}
 Let $g_1$ and $g_2$ be conformally equivalent metrics on
$M$. Then $(M,g_1)$ is conformally Osserman if and only if $(M,g_2)$ is conformally
Osserman.
\end{theorem}

Since $J_W(x)x=0$ and since $J_W(x)$ is self-adjoint, $J_W$ preserves the subspace $x^\perp$. We define the {\it
reduced conformal Jacobi operator} by letting
$$\tilde J_W(x):=J_W(x)|_{x^\perp}\,.$$
It is then immediate that $J_W$ has constant eigenvalues on $S_P(M,g)$ if and only if $\tilde J_W$ has constant
eigenvalues on $S_P(M,g)$; eliminating the trivial eigenvector $x$ simplifies subsequent statements.

Suppose that $(M,g)$ is locally conformally equivalent to a local 2 point homogeneous space $(M_0,g_0)$. Since the
local isometries of
$(M_0,g_0)$ act transitively on $S(M_0,g_0)$, the eigenvalues of $J_{W_0}$ are constant on $S(M_0,g_0)$
so $(M_0,g_0)$ is conformally Osserman. Thus by
Theorem \ref{thm-1.2},
$(M,g)$ is conformally Osserman. The following two results are partial converses to this observation.

We can classify conformally Osserman manifolds with certain eigenvalue structures:

\begin{theorem}\label{thm-1.3} Let $(M,g)$ be a conformally Osserman Riemannian manifold.
\begin{enumerate}
\item Suppose that $\tilde J_{W_P}$ has only one eigenvalue at each point $P$ of $M$. Then $(M,g)$ is conformally
flat.
\item Let $m\ge8$. Suppose that  $\tilde J_{W_P}$ has two distinct eigenvalues of multiplicities $1$ and $m-2$ for
each point
$P\in M$. Then
$(M,g)$ is locally conformally equivalent to either $(\mathbb{CP}^n,g_{FS})$ or
$({}^*\mathbb{CP}^n,{}^*g_{FS})$.
\end{enumerate}
\end{theorem}

We can use topological methods to control the eigenvalue structure in certain dimensions and derive the following
result from Theorem \ref{thm-1.3}:

\begin{theorem}\label{thm-1.4}
 Let $(M,g)$ be a conformally Osserman Riemannian manifold of dimension $m$.
\begin{enumerate}
\item If $m$ is odd, then $(M,g)$ is conformally flat.
\item If $m=4k+2\ge10$ and if $P$ is a point of $M$ where $W_P\ne0$, then there is an open neighborhood of $P$ in
$M$ which is conformally equivalent to an open subset of either
$(\mathbb{CP}^n,g_{FS})$ or $({}^*\mathbb{CP}^n,{}^*g_{FS})$.
\end{enumerate}
\end{theorem}

Here is a brief guide to the paper. In Section
\ref{sect-2} we prove Theorem \ref{thm-1.1} and classify the conformally complex space forms.
In Section \ref{sect-3}, we review results of Chi \cite{Ch88} in the algebraic
context. In Section
\ref{sect-4}, we establish Theorems
\ref{thm-1.3} and
\ref{thm-1.4} and thereby establish a conformal equivalent of the Osserman conjecture in certain situations.

\section{Conformally complex space forms}\label{sect-2}
We adopt an argument of Tricerri and Vanhecke \cite{TV81} to establish Theorem \ref{thm-1.1}. Let
$(M,g)$ be a conformally complex space form of dimension $m\ge8$. By assumption, there
exists a Hermitian almost complex structure $\Phi$ on $M$ so that
$$W=\lambda_0R_0+\lambda_1R_\Phi\quad\text{where}\quad\lambda_1\ne0\,.$$
If $g(x,x)=1$, then:
$$J_W(x)y=\left\{
\begin{array}{lll}
0&\text{if}&y=x.\\
(\lambda_0+3\lambda_1)y&\text{if}&y=\Phi x,\\
\lambda_0y&\text{if}&y\perp\{x,\Phi x\}\,.
\end{array}\right.
$$
Thus $\tilde J_W$ has two eigenvalues $(\lambda_0+3\lambda_1,\lambda_0)$ with multiplicities $(1,m-2)$. Since
$\Tr\{J_W\}=0$, this shows
\begin{eqnarray}
&&3\lambda_1+(m-1)\lambda_0=0\quad\text{so}\quad\lambda_0=-\ffrac3{m-1}\lambda_1\quad\text{and}\nonumber\\
&&R=\lambda_1R_\Phi+c_2(m)L+(c_1(m)\tau-\ffrac3{m-1}\lambda_1)R_0\,.\label{eqn-2.a}
\end{eqnarray}

If $g_1$ and $g_2$ are conformally related, then
$$
W_{g_1}=W_{g_2},\quad\Phi_{g_2}=\Phi_{g_1},\quad\text{and}\quad
\lambda_1(g_1)g_1=\lambda_1(g_2)g_2\,.
$$
Set $g_2:=|\lambda_1(g_1)|g_1$; $\lambda_1(g_2)=\pm1$. Thus we may therefore assume
henceforth without loss of generality that $\lambda_1=\pm1$. Let $\Phi_{;x}$, $W_{;x}$, and $R_{;x}$ denote the
covariant derivatives of these tensors. Since any complex space form is locally isometric to complex projective space with
a multiple of the Fubini-Study metric or to the negative curvature dual, Theorem
\ref{thm-1.1} will follow the following result:

\begin{lemma}\label{lem-2.1} Let $(M,g,\Phi)$ satisfy
$W=\lambda_1R_\Phi+\lambda_0R_0$
where $\lambda_1=\pm1$. Assume that $m\ge8$. Let $\{a,\Phi a,b,\Phi b,c,\Phi c\}$ be an
orthonormal set. Then:\begin{enumerate}
\item We have $\Phi_{;a}\Phi=-\Phi\Phi_{;a}$. We also have $\Phi$, $\Phi_{;a}$, and $\Phi\Phi_{;a}$ are
skew-adjoint.
\item We have $(\Phi_{;c}b-\Phi_{;b}c,a)=0$.
\item We have $\Phi_{;a}a=0$,
$\Phi_{;a}\Phi a=0$, and $\Phi_{;a}b+\Phi_{;b}a=0$.
\item We have $\nabla\Phi=0$.
\item We have that $(M,g,\Phi)$ is a complex space form.
\end{enumerate}
\end{lemma}

\begin{proof} We covariantly differentiate the identity $\Phi^2=-\id$ to see $\Phi_{;a}\Phi+\Phi\Phi_{;a}=0$. As
$\Phi$ is skew-adjoint, $\Phi_{;a}$ is skew-adjoint. The fact that $\Phi\Phi_{;a}$ is skew-adjoint
then follows from the fact that $\Phi_{;a}$ and $\Phi$ anti-commute. Assertion (1) follows.

We use the second Bianchi identity
\begin{equation}\label{eqn-2.b}
\{R_{;x}(y,z)+R_{;y}(z,x)+R_{;z}(x,y)\}w=0
\end{equation}
to prove Assertions (2) and (3). Let $\sigma_{x,y,z}$ be summation with
respect to the cyclic permutation of $(x,y,z)$. Equations (\ref{eqn-2.a}) and (\ref{eqn-2.b}) imply:
\begin{equation}\label{eqn-2.c}
\begin{array}{l}
0=\sigma_{x,y,z}\{c_1\tau_{;x}[g(z,w)y-g(y,w)z]\\
\quad+ c_2 [g(\rho_{;x}z,w)y-g(\rho_{;x}y,w)z +g(z,w)\rho_{;x}y-g(y,w)\rho_{;x}z]\vphantom{\vrule height
11pt}\\
\quad+\lambda_1[g(\Phi_{;x}z,w)\Phi y+g(\Phi z,w)\Phi_{;x}
y-g(\Phi_{;x}y,w)\Phi z-
                 g (\Phi y,w)\Phi_{;x} z\vphantom{\vrule height 11pt}\\
\quad-2g(\Phi_{;x}y,z)\Phi w-2g(\Phi y,z)\Phi_{;x} w]\}\,.\vphantom{\vrule height 11pt}
\end{array}\end{equation}

Since $m\geq 8$, we may choose  $d$ so that $\{a,\Phi a,b,\Phi b,c,\Phi
c,d,\Phi d\}$ is an orthonormal set. Let $x=a$, $y=b$, $z=c$, and $w=d$ in Equation (\ref{eqn-2.c}). Then:
\begin{eqnarray*}
0&=&\sigma_{a,b,c}\{c_2 [g(\rho_{;a}c,d)b-g(\rho_{;a}b,d)c]\\
 &+&\lambda_1[g(\Phi_{;a}c,d)\Phi b-g(\Phi_{;a}b,d)\Phi c-2g(\Phi_{;a}b,c)\Phi d]\}\,.
\end{eqnarray*}
Since $\lambda_1\ne0$, $g(\Phi_{;a}c-\Phi_{;c}a,d)\Phi b=0$.
Assertion (2) follows by setting the coefficient of $\Phi b$ to
zero and permuting $a$, $b$, $c$ and  $d$ appropriately.

To prove Assertion (3), set $x=a$, $y=b$, $z=\Phi b$, and $w=d$ in Equation (\ref{eqn-2.c}):
\begin{eqnarray*}
0&=&c_1\tau_{;a}[g(\Phi b,d)b-g(b,d)\Phi b]\\
&+& c_2 [g(\rho_{;a}\Phi b,d)b-g(\rho_{;a}b,d)\Phi b +g(\Phi b,d)\rho_{;a}b-g(b,d)\rho_{;a}\Phi b]\\
&+&\lambda_1[g(\Phi_{;a}\Phi b,d)\Phi b+g(\Phi \Phi b,d)\Phi_{;a}
b-g(\Phi_{;a}b,d)\Phi \Phi b-
                 g (\Phi b,d)\Phi_{;a} \Phi b\\
&-&2g(\Phi_{;a}b,\Phi b)\Phi d-2g(\Phi b,\Phi b)\Phi_{;a} d]\\
&+&c_1\tau_{;b}[g(a,d)\Phi b-g(\Phi b,d)a]\\
&+& c_2 [g(\rho_{;b}a,d)\Phi b-g(\rho_{;b}\Phi b,d)a +g(a,d)\rho_{;b}\Phi b-g(\Phi b,d)\rho_{;b}a]\\
&+&\lambda_1[g(\Phi_{;b}a,d)\Phi \Phi b+g(\Phi a,d)\Phi_{;b} \Phi
b-g(\Phi_{;b}\Phi b,d)\Phi a-
                 g (\Phi \Phi b,d)\Phi_{;b} a\\
&-&2g(\Phi_{;b}\Phi b,a)\Phi d-2g(\Phi \Phi b,a)\Phi_{;b} d]\\
&+&c_1\tau_{;\Phi b}[g(b,d)a-g(a,d)b]\\
&+& c_2 [g(\rho_{;\Phi b}b,d)a-g(\rho_{;\Phi b}a,d)b +g(b,d)\rho_{;\Phi b}a-g(a,d)\rho_{;\Phi b}b]\\
&+&\lambda_1[g(\Phi_{;\Phi b}b,d)\Phi a+g(\Phi b,d)\Phi_{;\Phi b}
a-g(\Phi_{;\Phi b}a,d)\Phi b-
                 g (\Phi a,d)\Phi_{;\Phi b} b\\
&-&2g(\Phi_{;\Phi b}a,b)\Phi d-2g(\Phi a,b)\Phi_{;\Phi b} d]\,.
\end{eqnarray*}
which simplifies to become:
\begin{eqnarray}
0&=&c_2 [g(\rho_{;a}\Phi b,d)b-g(\rho_{;a}b,d)\Phi
b]+c_2[g(\rho_{;b}a,d)\Phi b-g(\rho_{;b}\Phi b,d)a]
   \nonumber\\
&+&c_2[g(\rho_{;\Phi b}b,d)a-g(\rho_{;\Phi b}a,d)b]\nonumber\\
&+&\lambda_1[g(\Phi_{;a}\Phi b,d)\Phi
b+g(\Phi_{;a}b,d)b-2g(\Phi_{;a}b,\Phi b)\Phi d-2\Phi_{;a}d]\}
   \label{eqn-2.d}\\
&+&\lambda_1[-g(\Phi_{;b}a,d)b-g(\Phi_{;b}\Phi b,d)\Phi
a-2g(\Phi_{;b}\Phi b,a)\Phi d]\nonumber\\
&+&\lambda_1[(g(\Phi_{;\Phi b}b,d)\Phi a-g(\Phi_{;\Phi b}a,d)\Phi
b-2g(\Phi_{;\Phi b}a,b)\Phi d]\,.\nonumber
\end{eqnarray}

We apply Assertions (1) and (2) to see
\begin{eqnarray*}
0&=&\{g(\rho_{;a}\Phi b,d)-g(\rho_{;\Phi b}a,d)\}b
  =\{-g(\rho_{;a}b,d)+g(\rho_{;b}a,d)\}\Phi b,\\
0&=&\{g(\Phi_{;a}\Phi b,d)-g(\Phi_{;\Phi b}a,d)\}\Phi b
  =\{g(\Phi_{;a}b,d)-g(\Phi_{;b}a,d)\}b\\
0&=&g(\Phi_{;a}b,\Phi b),\quad\text{and}\quad g(\Phi_{;\Phi b}a,b)=-g(\Phi_{;\Phi b}b,a)\,.
\end{eqnarray*}
Consequently, we may rewrite Equation (\ref{eqn-2.d}) in the form:
\begin{eqnarray*}
0&=&c_2[-g(\rho_{;b}\Phi b,d)a+g(\rho_{;\Phi b}b,d)a]-2\lambda_1\Phi_{;a}d\\
&-&2\lambda_1\{g(\Phi_{;b}\Phi b,a)-g(\Phi_{;\Phi b}b,a)\}\Phi d\\
&+&\lambda_1\{-g(\Phi_{;b}\Phi b,d)+g(\Phi_{;\Phi b}b,d)\}\Phi a\,.
\end{eqnarray*}
As by Assertion (1) $\Phi\Phi_{;a}$ is skew-adjoint, $g(\Phi_{;a}d,\Phi d)=0$. Taking the inner
products with
$\Phi d$ and with
$\Phi a$ then yields
\begin{eqnarray}
&&0=g(\Phi_{;b}\Phi b-\Phi_{;\Phi b}b,a),\label{eqn-2.e}\\
&&0=g(\Phi_{;\Phi b}b-\Phi_{;b}\Phi
b,d)-2\lambda_1g(\Phi_{;a}d,\Phi a)\,.\label{eqn-2.f}
\end{eqnarray}
Equation (\ref{eqn-2.e}) shows that $(\Phi_{;b}\Phi b-\Phi_{;\Phi b}b)\in\operatorname{Span}\{b,\Phi b\}$. Since
\begin{equation}\label{eqn-2.g}
g(\Phi_{\Phi b}b,b)=g(\Phi_{\Phi b}b,\Phi b)=g(\Phi_{;b}\Phi b,b)=g(\Phi_{;b}\Phi b,\Phi b)=0,
\end{equation}
we may conclude $\Phi_{\Phi b}b-\Phi_{;b}\Phi b=0$. Equation (\ref{eqn-2.f}) then
implies
$$0=-2g(\Phi_{;a}d,\Phi a)=2g(d,\Phi_{;a}\Phi a)\,.$$
This implies $\Phi_{;a}\Phi a\in\operatorname{Span}\{a,\Phi a\}$ and applying Equation (\ref{eqn-2.g}) shows
$\Phi_{;a}\Phi a=0$. We then see $\Phi_{;a}a=-\Phi\Phi\Phi_{;a}a=\Phi\Phi_{;a}\Phi a=\Phi
0=0$. The final identity of Assertion (3) then follows by polarization.

By Assertions (2) and (3),  $g(\Phi_{;a}b,c)=0$ so $\Phi_{;a}b\in\operatorname{Span}\{a,\Phi a,b,\Phi b\}$. We
show $\Phi_{;a}b=0$ and establish Assertion (4) by computing:
\begin{eqnarray*}
&&0=g(\Phi_{;a}b,a)=-g(\Phi_{;a}a,b)=0,\\
&&0=g(\Phi_{;a}b,\Phi a)=-g(\Phi_{;a}\Phi a,b)=0,\\
&&0=g(\Phi_{;a}b,b)=g(\Phi_{;a}b,\Phi b)\,.
\end{eqnarray*}

Because $\nabla\Phi=0$, $R(x,y)\Phi=\Phi R(x,y)$. We compute:
\begin{equation}\label{eqn-2.h}\begin{array}{l}
R(x,y)\Phi z=\lambda_1\{g(y,z)\Phi x-g(x,z)\Phi y+2g(\Phi x,y)z\}\\
\quad+c_2(m)\{g(\rho y,\Phi z)x-g(\rho x, \Phi z)y+g(y,\Phi z)\rho
x-g(x,\Phi z)\rho y\}
   \vphantom{\vrule height11pt}\\
\quad+(c_1(m)\tau-\ffrac3{m-1}\lambda_1)\{g(y,\Phi z)x-g(x,\Phi z)y\}
   \vphantom{\vrule height11pt}\\
\Phi R(x,y)z=\lambda_1\{-g(\Phi y,z)x+g(\Phi x,z)y+2g(\Phi x,y)z\}
   \vphantom{\vrule height11pt}\\
\quad+c_2(m)\{g(\rho y,z)\Phi x-g(\rho x,z)\Phi y+g(y,z)\Phi\rho
x-g(x,z)\Phi\rho y\}
   \vphantom{\vrule height11pt}\\
\quad+(c_1(m)\tau-\ffrac3{m-1}\lambda_1)\{g(y,z)\Phi x-g(x,z)\Phi y\}\,.
   \vphantom{\vrule height11pt}
\end{array}\end{equation}
Let $x=a$, $y=b$, and $z=a$ in Equation (\ref{eqn-2.h}):
\begin{eqnarray}
&&\lambda_1\{-\Phi b\}+c_2(m)\{g(\rho b,\Phi a)a-g(\rho a,\Phi a)b \}\nonumber\\
&=&c_2(m)\{g(\rho b,a)\Phi a-g(\rho a,a)\Phi b-\Phi\rho
b\}+(c_1(m)\tau-\ffrac3{m-1}\lambda_1)\{-\Phi
b\}\,.\label{eqn-2.i}
\end{eqnarray}
This implies $\Phi\rho b\in\operatorname{Span}\{\Phi b,\Phi
a,b,a\}$. Similarly $\Phi\rho b\in\operatorname{Span}\{\Phi b,\Phi
c,b,c\}$. This implies $\rho b\in\operatorname{Span}\{b,\Phi b\}$.
Taking the inner product of Equation (\ref{eqn-2.i}) with $b$ we get
$-g(\rho a, \Phi a)=-g(\Phi\rho b,b)=g(\rho b,\Phi b)$. Symmetrizing over $\{a,b,c\}$, we see $g(\rho b,\Phi b)=0$ so
$\rho b\in\operatorname{Span}\{b\}$ and hence $\rho b=\mu(b)b$. As
\begin{eqnarray*}
&&\rho b=\mu(b)b,\quad \rho a=\mu(a)a,\quad\text{and},\\
&&\rho(a+b)=\mu(a+b)(a+b)=\mu(a)a+\mu(b)b,
\end{eqnarray*}
$\mu(a+b)=\mu(a)=\mu(b)$ and hence $\mu$ is constant. This implies $(M,g)$ is Einstein.
Consequently, $L_0$ is a multiple of $R_0$.
 Equation (\ref{eqn-2.a}) now implies that $(M,g,\Phi)$ is a complex space form.
\end{proof}

\section{Algebraic curvature tensors}\label{sect-3}

We consider a triple $\mathcal{V}:=(V,g,A)$ where $g$ is a positive definite inner
product on a  real vector space $V$ of dimension $m$ and where
$A\in\otimes^4V^*$ is an {\it algebraic curvature tensor} on $V$; i.e. $A$ has the usual symmetries of the Riemann
curvature tensor:
\begin{eqnarray*}
&&A(x,y,z,w)=A(z,w,x,y)=-A(y,x,z,w),\quad\text{and}\\
&&A(x,y,z,w)+A(y,z,x,w)+A(z,x,y,w)=0\,.
\end{eqnarray*}

We follow the discussion in \cite{Gi02}. If $\Psi$ is a self-adjoint map of $V$, set:
\begin{equation}\label{eqn-3.a}
A_\Psi(x,y,z,w):=g(\Psi x,w)g(\Psi y,z)-g(\Psi x,z)g(\Psi y,w)\,.
\end{equation}
For example, the algebraic curvature tensors $L$ and $R_0$ of Equation (\ref{eqn-1.a}) can be expressed in the
form:
$$R_0=R_{\id}\quad\text{and}\quad L=R_{\id+\rho}-R_{\id}-R_\rho\,.$$
 Similarly if $\Phi$ is skew-adjoint, we generalize Equation (\ref{eqn-1.d}) and set:
\begin{equation}\label{eqn-3.b}
A_\Phi(x,y,z,w):=g(\Phi x,w)g(\Phi y,z)-g(\Phi x,z)g(\Phi y,w)-2g(\Phi x,y)g(\Phi z,w)\,.
\end{equation}
One checks easily that the tensors $A_\Psi$ and $A_\Phi$ defined above are algebraic curvature tensors. One has
the following result of Fiedler \cite{F01}:

\begin{theorem}\label{thm-3.1}
 The space of all algebraic curvature tensors is a real vector space. It is
spanned by the tensors $A_\Psi$ of Equation (\ref{eqn-3.a}). It is also spanned by the tensors
$A_\Phi$ of Equation (\ref{eqn-3.b}).
\end{theorem}

If $A$ is an algebraic curvature tensor, then we define the associated curvature operator $A(x,y)$ and Jacobi
operator
$J_A(x)$ by the relations:
$$g(A(x,y)z,w)=A(x,y,z,w)\quad\text{and}\quad g(J_A(x)y,z)=A(y,x,x,z)\,.$$
The operator $J_A(x)=A(\cdot,x)x$ is a self-adjoint map of $V$ and we say that $\mathcal{V}$ is {\it Osserman} if
the eigenvalues of
$J_A$ are constant on $S(V):=\{x\in V:g(x,x)=1\}$. We note that $J_A(x)x=0$ and let $\tilde J_A(x)$ be
the restriction of $J_A(x)$ to $x^\perp$.

The following classification result is due to Chi
\cite{Ch88} and will be crucial in establishing Theorem \ref{thm-1.3}.
\begin{theorem}\label{thm-3.2}
 Let $A$ be an Osserman algebraic curvature tensor on $\mathcal{V}$.
\begin{enumerate}
\item If $\tilde J_A$ has only one eigenvalue, then $A=\lambda A_{\id}$.
\item If $\tilde J_A$ has two eigenvalues and if one of those eigenvalues has multiplicity $1$, then there exists
a Hermitian almost complex structure $\Phi$ on $V$ and there exist real constants $\lambda_0$ and $\lambda_1$ so
that
$A=\lambda_0A_{\id}+\lambda_1A_\Phi$.
\end{enumerate}
\end{theorem}

The following observation is also due to Chi \cite{Ch88}; it is a straightforward application of work of Adams
\cite{Adams62} concerning vector fields on spheres and will be critical in proving Theorem \ref{thm-1.4}:

\begin{theorem}\label{thm-3.3}
 Let $A$ be an Osserman algebraic curvature tensor on $\mathcal{V}$.
\begin{enumerate}
\item If $m$ is odd, then $\tilde J_A$ has only one eigenvalue.
\item If $m\equiv2$ mod $4$, then either $\tilde J_A$ has only eigenvalue or $\tilde J_A$ has exactly $2$
eigenvalues and one of those eigenvalues has multiplicity $1$.
\end{enumerate}
\end{theorem}

\section{Conformal Osserman manifolds}\label{sect-4}

\begin{proof}[Proof of Theorem \ref{thm-1.3} (1)] Assume that $(M,g)$ is conformally Osserman and that $\tilde J_W$
has only one eigenvalue. By Theorem \ref{thm-3.2}, $W_P=\lambda_PA_{\id}$. By Equation
(\ref{eqn-1.e}),
$$0=\Tr\{J_W(x)\}=(m-1)\lambda_Pg(x,x)\,.$$
Thus $\lambda_P=0$ so $W=0$ and $(M,g)$ is conformally flat.
\end{proof}

\begin{proof}[Proof of Theorem \ref{thm-1.3} (2)] Assume that
$(M,g)$ is conformally Osserman,  that $\tilde J_W$ has two distinct eigenvalues,
that one of the eigenvalues has multiplicity $1$, and that $m\ge8$. By Theorem \ref{thm-3.2}, there
exists a Hermitian almost complex structure $\Phi$ on each tangent space so that
$$W=\lambda_0R_{\id}+\lambda_1R_\Phi\quad\text{where}\quad\lambda_0=-\ffrac3{m-1}\lambda_1\,.$$
We then use techniques developed in \cite{Gi02} to show that $\Phi$ can be chosen to vary smoothly with $P$, at
least locally. Thus
$(M,g)$ is a conformally complex space form so Theorem
\ref{thm-1.3} (2) follows from Theorem \ref{thm-1.1}.
\end{proof}

\begin{proof}[Proof of Theorem \ref{thm-1.4} (1)] Let  $(M,g)$ be conformally Osserman. If $m$ is odd, then
Theorem
\ref{thm-3.3} implies that $\tilde J_W$ has only one eigenvalue. Therefore by Theorem
\ref{thm-1.3}, $(M,g)$ is conformally flat.\end{proof}

\begin{proof}[Proof of Theorem \ref{thm-1.4} (2)] Let $m=4k+2\ge10$ and let $P$ be a point of $M$ where $W_P\ne0$.
By Theorem
\ref{thm-3.3}, either
$\tilde J_{W_P}$ has only one eigenvalue or $\tilde J_{W_P}$ has two
eigenvalues and one has multiplicity $1$. If $\tilde J_{W_P}$ has only one eigenvalue, then that eigenvalue is $0$
by Equation (\ref{eqn-1.e}). This implies $J_{W_P}=0$ and hence $W_P=0$ which is contrary to the assumption which
we have made. Consequently $\tilde J_{W_P}$ has two eigenvalues at $P$ and hence on a neighborhood $\mathcal{O}_P$.
By Theorem \ref{thm-1.3} (2), $(\mathcal{O}_P,g)$ is locally conformally equivalent to an open subset of either
complex projective space with the Fubini study metric or to the negative curvature dual.\end{proof}

\section*{Acknowledgments} Research of N. Bla{\v z}i{\'c} partially supported by the DAAD (Germany)
 and MNTS Project \#1854 (Srbija). Research of
P. Gilkey partially supported by the MPI (Leipzig).
We thank  S. Nik{\v c}evi{\'c} for her interest for this work and useful comments.
The authors wish to express their thanks to the Technical University of Berlin where much of the research
reported here was conducted.

\end{document}